\documentclass[12pt,reqno,oneside]{amsart}

\newtheorem{defin}{Definition}[section]
\newtheorem{theorem}[defin]{Theorem} 
\newtheorem{corollary}[defin]{Corollary}
\newtheorem{lemma}[defin]{Lemma} 

\def\C{{\mathbb  C}}
\def\N{{\mathbb  N}}
\def\R{{\mathbb  R}}
\def\cR{{\mathcal R}}

\def\d{\delta}
\def\e{\epsilon}
\def\z{\zeta}

\def\l{\lambda}

\def\bar{\overline}        

\def\disc{\triangle}             
\def\cdisc{\bar\triangle}        
\def\bdisc{b\disc}               

\def\ctwo{{\mathbb C}^2}

\def\st{such that}
\def\smnt{Stein manifold}
\def\psh{plurisubharmonic}               
\def\spsh{strongly\ plurisubharmonic}
\def\cont{continuous}
\def\holo{holomorphic}

\def\ss{\subset\!\subset}   

\begin{document}
\begin{flushright}
\end{flushright}

\medskip

\title{\bf Proper discs in Stein manifolds avoiding complete pluripolar sets}         

\author{Barbara Drinovec Drnov\v sek}

\bigskip\rm

\maketitle
\pagestyle{plain}
\section{Introduction and the results}
\baselineskip 16pt
Denote by $\disc$ the open unit disc in $\C$.
Recall that a subset $Y$ in a complex manifold $X$ is called {\it complete pluripolar} if
there exists a \psh\ function $\rho$ on $X$ \st\ $Y=\{z;\,\rho(z)=-\infty\}$.

In this paper we prove the following result.

\begin{theorem}
Let $X$ be a \smnt\ of dimension at least $2$. Given a 
closed complete pluripolar set $Y\subset X$, a point $p\in X\setminus Y$ 
and a vector $v$ tangent to $X$ at $p$, 
there exists a proper holomorphic map $f\colon \disc\to X$ \st\ $f(0)=p$,
$f'(0)=\l v$ for some $\l>0$ 
and $f(\disc)\cap Y=\emptyset$.
\label{thm1.1}
\end{theorem}

Clearly, every closed complex analytic subset $A$ of a connected 
Stein manifold $X$, $A\neq X$, is locally complete pluripolar, that is,
for any point $a\in A$ there is an open neighborhood $U$ of $a$ such that
$A\cap U$ is complete pluripolar in $U$. By \cite{Col} every closed locally complete
pluripolar set in a Stein manifold is complete pluripolar, 
thus every closed complex analytic subset
is closed complete pluripolar.
Therefore our theorem answers the question posed in \cite{FG} on the
existence of proper holomorphic discs in the complements of hypersurfaces.

J. Globevnik \cite{Glo} proved in 2000 that for any point $p$ in a \smnt\ $X$
of dimension at least $2$ there exists a proper holomorphic map from the
unit disc to $X$ with the point  $p$ in its image.

The most general result on avoiding certain sets by proper holomorphic discs was 
given by H. Alexander \cite{Ale} in 1975: he proved that for a closed
polar set $E\subset \C$ there exists a proper holomorphic map 
$F=(F_1,F_2)\colon\disc\to\ctwo$ \st\ $F_1(\disc)\cap E=\emptyset$. 
On the other hand, 
a proper holomorphic disc in $\ctwo$ cannot avoid a non-polar 
set of parallel complex lines (see \cite{Jul,Tsu,Ale,FG}). 
F. Forstneri\v c and J. Globevnik \cite{FG} in 2001 constructed a proper holomorphic
disc in $\ctwo$ omitting both coordinate axes and proper holomorphic 
discs avoiding large real cones in $\ctwo$.
However, it was unknown if the image of a proper holomorphic map from the disc 
can miss three or more complex lines.
Our theorem provides a positive answer to this question since
a finite union of complex lines in $\ctwo$ is closed
complete pluripolar.
Note that
closed convex sets in $\ctwo$ which can be avoided by the image of proper holomorphic 
maps from the disc were characterized in \cite{Dri}.

We shall prove the following approximation theorem, which easily implies Theorem
\ref{thm1.1}. In fact, Theorem \ref{thm1.1} will follow directly from Lemma
\ref{lema}.

\begin{theorem}
Let $X$ be a \smnt\ of dimension at least $2$ and let $Y\subset X$
be a closed complete pluripolar subset. Let $d$ be a complete metric on $X$ which induces 
the manifold topology.
Assume that $f\colon\disc\to X$ 
is a holomorphic map 
\st\ there is an open subset $V\ss\disc$ with the property
$f(\z)\notin Y$ for $\z\in\disc\setminus V$.
Given $\epsilon>0$ there is
a proper holomorphic map $g\colon\disc\to X$ satisfying
\item[(i)] $g(\z)\notin Y$ for $\z\in\disc\setminus V$,
\item[(ii)] $d(g(\z),f(\z))<\e$ for $\z\in V$,
\item[(iii)] $g(0)=f(0)$ and $g'(0)=\lambda f'(0)$ for some $\l>0$.
\label{thm1.2}
\end{theorem}
 
We will prove Theorem \ref{thm1.2} in section 2. 
 
\begin{corollary}
Let $X$ be a \smnt\ of dimension at least $2$ and let $Y\subset X$
be a closed complete pluripolar subset. Assume that $S$ is a discrete subset of $X$ 
\st\ $S\cap Y=\emptyset$. Then there are proper holomorphic maps $f_n\colon\disc\to X$
\st\ $f_n(\disc)$ are pairwise disjoint, $f_n(\disc)$ avoids $Y$ $(n\in\N)$
and $\cup_nf_n(0)=S$.
\label{cor1}
\end{corollary}

\begin{proof}
We first note that a finite union of complete pluripolar sets is complete pluripolar,
since a finite sum of \psh\ functions is \psh.
We will also need the fact that a discrete set $S$ in a Stein manifold is complete pluripolar.
Namely, by \cite{Col} it is enough to prove that $S$ is locally complete pluripolar,
which follows from the fact that $S$ is a complex analytic subset of $X$.

Let $S=\{s_n;\, n\in\N\}$. We shall construct the maps $f_n$ inductively.
By Theorem \ref{thm1.1} there is a proper holomorphic map 
$f_1\colon\disc\to X$
\st\ $f_1(\disc)\cap (Y\cup S\setminus\{s_1\})=\emptyset$ 
and $f_1(0)=s_1$.
Assume that for some $n\in\N$ we have already constructed proper holomorphic maps
$f_j\colon\disc\to X$, $1\le j\le n$,  such that
$f_j(\disc)$ are pairwise disjoint, $f_j(\disc)$ avoids $Y$
and $f_j(\disc)\cap S=\{s_j\}$  $(1\le j\le n)$.
By Remmert's proper mapping theorem \cite{Re1,Re2}, \cite[p. 65]{Ch2}
the image of a proper holomorphic map is a closed analytic 
subset of $X$ and therefore closed complete pluripolar. 
Thus
$Y\cup f_1(\disc)\cup\cdots\cup f_n(\disc)\cup S\setminus\{s_{n+1}\}$ is
closed complete pluripolar. Then by Theorem \ref{thm1.1} there is a proper holomorphic
map $f_{n+1}\colon\disc\to X$ \st\ 
$$f_{n+1}(\disc)\cap (Y\cup f_1(\disc)\cup\cdots\cup f_n(\disc)\cup S\setminus\{s_{n+1}\})=
\emptyset$$ and
$f_{n+1}(0)=s_{n+1}$.
The inductive construction is finished and the proof is complete.
\end{proof}

Let $\cR$ be a bordered Riemann surface. 
By the theorem of Ahlfors \cite{Ahl}, there are inner functions on $\cR$.
Recall that a nonconstant continuous function $f\colon\cR\to\cdisc$,
which is holomorphic on $\cR\setminus b\cR$,
is called an {\it inner function }
(or an {\it Ahlfors function}) on $\cR$ 
if $|f|=1$ on $b\cR$.  Therefore Theorem \ref{thm1.1} implies the following:

\begin{corollary}
Let $X$ be a \smnt\ of dimension at least $2$ and let $Y\subset X$
be a closed complete pluripolar subset. 
Given a bordered Riemann surface $\cR$ there is a proper holomorphic
map $f\colon\cR\setminus b\cR\to X$ \st\ $f(\cR\setminus b\cR)\cap Y=\emptyset$.
\end{corollary}

\section{Proof of Theorem 1.2}

As it was observed in \cite{FG} the methods developed in \cite{FG1,Glo}  actually
prove the following:
\begin{theorem}
Let $X$ be a \smnt\ of dimension at least $2$ and $\rho\colon X\to\R$
a smooth exhaustion function which is \spsh\ on $\{\rho>M\}$ for some $M\in\R$.
Let $f\colon\cdisc\to X$ be a continuous map which is holomorphic on $\disc$
\st\ $\rho(f(\z))>M$ for each $\z\in\bdisc$. Let $d$ be a complete metric on $X$
which induces the manifold topology.
For any numbers $0<r<1$, $\e>0$, $N>M$ and for any finite set $A\subset \disc$
there exists a \cont\ map $g\colon\cdisc\to X$, \holo\ on $\disc$, satisfying
\item[(i)] $\rho(g(\z))>N$ for $\z\in\bdisc$,
\item[(ii)] $\rho(g(\z))>\rho(f(\z))-\e$ for $\z\in\cdisc$,
\item[(iii)] $d(f(\z),g(\z))<\e$ for $|\z|\le r$, and
\item[(iv)] $g(\z)=f(\z)$ and $g'(\z)=f'(\z)$ for $\z\in A$.
\label{thmFG}
\end{theorem}

In the proof of Theorem \ref{thm1.2} we also need the following lemma
which is a slight generalization of \cite[Lemma 1]{Ch1}. Since its proof
is essentially the same, we omit it.

\begin{lemma}
Let $X$ be a \smnt\ and $Y\subset X$ a 
complete pluripolar set.
Let $L_1\subset L_2\subset X$ be holomorphically convex compact sets.
Then the set $(L_1\cup Y)\cap L_2$ is holomorphically convex.
\label{lema2.3}
\end{lemma}

The main tool in the proof of Theorem \ref{thm1.2} is the following

\begin{lemma}
Let $X$ be a \smnt\ of dimension at least $2$ and let $Y\subset X$
be a closed complete pluripolar subset. Let $d$ be a complete metric on $X$ which induces 
the manifold topology.
Assume that $f\colon\disc\to X$ is a holomorphic map 
\st\ there is an open subset $V\ss\disc$ with the property
$f(\z)\notin Y$ for $\z\in\disc\setminus V$.
Given $\epsilon>0$ there are a domain $\Omega$, $\{0\}\cup V\ss\Omega\ss\disc$,
conformally equivalent to the unit disc and
a proper holomorphic map $g\colon\Omega\to X$ with the following properties
\item[(i)] $g(\z)\notin Y$ for $\z\in\Omega\setminus V$,
\item[(ii)] $d(g(\z),f(\z))<\e$ for $\z\in V$,
\item[(iii)] $g(0)=f(0)$ and $g'(0)=f'(0)$.
\label{lema}
\end{lemma}

\begin{proof}
One can choose a simply connected domain $\Omega_1$ \st\ $\{0\}\cup V\ss\Omega_1\ss\disc$.
By \cite[Theorem 5.1.6]{Hor} there is a smooth \spsh\ exhaustion function $\rho$ 
for Stein manifold $X$.
Sard's theorem implies that one can choose a strictly increasing sequence $\{M_n\}$ 
of regular values of $\rho$ converging to $\infty$ with
$M_1$ so big that $\rho(f(\z))<M_1$ for $\z\in \bar\Omega_1$.
By continuity 
there is a simply connected domain $\disc_1$, $\Omega_1\ss\disc_1\ss\disc$, \st\
$\rho(f(\z))<M_1$ for $\z\in \cdisc_1$.
Let $U_0=\emptyset$ and for $n\in\N$ 
denote by $U_n$ the sublevel set $\{z\in X;\, \rho(z)<M_n\}$.
Since $M_n$ is a regular value of $\rho$ it holds that 
$\bar U_n=\{z\in X;\, \rho(z)\le M_n\}$. This implies that 
$\bar U_n$ is a holomorphically convex compact set, because on a Stein manifold
plurisubharmonic hull equals holomorphic hull.

We shall construct inductively
a decreasing sequence of domains $\{\disc_n\}$ conformally 
     equivalent to $\disc$, an increasing sequence 
     of domains $\{\Omega_n\}$ conformally 
     equivalent to $\disc$,
$\Omega_n\ss\disc_n$ $(n\in\N)$,     
a sequence of \cont\ maps $g_n\colon\cdisc_n\to X$,
     \holo\ on $\disc_n$, and
a decreasing sequence of positive numbers $\{\e_n\}$,
satisfying for each $n\in\N$ the following
\begin{enumerate}
\item[(I)]  $g_n(\z)\in U_{n}\setminus(\bar U_{n-1}\cup Y)$\ \  
$(\z\in\bar{\disc_n\setminus \Omega_n})$,
\item[(II)] $g_{n+1}(\z)\notin \bar U_{n-1}\cup Y$\ \  $(\z\in\bar{\disc_{n+1}\setminus \Omega_n})$,
\item[(III)]$d(g_{n}(\z),g_{n+1}(\z))<\tfrac{\e_n}{2^n}$\ \  $(\z\in \Omega_{n})$,
\item[(IV)] $g_{n+1}(0)=g_{n}(0)$ and $g_{n+1}'(0)=g_{n}'(0)$,
\item[(V)]  if $z\in X$ \st\ $d(z,g_n(\bar{\disc_{n}\setminus V}))<\e_n$ then $z\in U_n\setminus Y$,
\item[(VI)] if $z\in X$ \st\ $d(z,g_{n+1}(\bar{\disc_{n+1}\setminus \Omega_n}))<\e_{n+1}$
            then $z\notin \bar U_{n-1}$.
\end{enumerate}

Let $g_1=f$ and let $\disc_1$ and $\Omega_1$ as above. 
Then (I) holds. Choose $\e_1$, $0<\e_1<{\e}$, so small that
(V) holds for $n=1$.
Suppose that $j\in \N$ and that we have constructed
$g_n$, $\disc_n$, $\Omega_n$ and $\e_n$, $1\le n\le j$, \st\
(I) and (V) hold for $1\le n\le j$ and (II), (III), (IV) and (VI)
hold for $1\le n\le j-1$.
It follows by Lemma \ref{lema2.3} that the set $(\bar U_{j-1}\cup Y)\cap \bar U_{j+1}$ is
holomorphically convex. Therefore there is a 
smooth \spsh\ exhaustion function $\rho_{j+1}$ on $X$ \st\ 
$\rho_{j+1}(z)<0$ $(z\in (\bar U_{j-1}\cup Y)\cap \bar U_{j+1})$  
and $\rho_{j+1}(g_j(\z))>0$ $(\z\in\bar{\disc_j\setminus \Omega_j})$ \cite[Theorem 5.1.6]{Hor}.
There is $N$ so big that $U_{j}\subset \{z;\,\rho_{j+1}(z)<N\}$.
We use Theorem \ref{thmFG} to get a \cont\ map
$g_{j+1}\colon\cdisc_j\to X$, holomorphic on $\disc_j$, with the following properties
\begin{enumerate}
\item[(a)] $\rho_{j+1}(g_{j+1}(\z))>N$ for $\z\in\bdisc_j$,
\item[(b)] $\rho_{j+1}(g_{j+1}(\z))>0$ for $\z\in\bar{ \disc_j\setminus \Omega_{j}}$,
\item[(c)] $d(g_{j+1}(\z),g_{j}(\z))<\tfrac{\e_j}{2^j}$ for $\z\in \Omega_{j}$, and
\item[(d)] $g_{j+1}(0)=g_{j}(0)$ and $g_{j+1}'(0)=g_{j}'(0)$.
\end{enumerate}

By (a) and by the choice of $N$ we get that $\rho(g_{j+1}(\z))>M_{j}$ $(\z\in\bdisc_{j})$.
Thus there is $M$, $M_{j}<M<M_{j+1}$, \st\ the holomorphic disc $g_{j+1}(\disc_j)$ and the level set 
$\{z;\,\rho(z)=M\}$ intersect transversally. 
It follows by (c) and (V) that $g_{j+1}(\Omega_{j}\setminus V)\subset U_j$.
This and the fact that $\rho\circ g_{j+1}$ is subharmonic imply that
there is a simply connected component of
$\{\z\in\disc_j;\, \rho(g_{j+1}(\z))<M\}$ which contains $\Omega_{j}$. We denote such 
component by $\disc_{j+1}$.
Choose a simply connected domain $\Omega_{j+1}$, $\Omega_j\ss \Omega_{j+1}\ss\disc_{j+1}$, 
\st\ $\rho(g_{j+1}(\z))>M_j$ $(\z\in\bar {\disc_{j+1}\setminus \Omega_{j+1}})$.
It is easy to see that 
$\disc_{j+1}$, $\Omega_{j+1}$, $g_{j+1}$ satisfy (I) for $n=j+1$ and
(II), (III) and (IV) for $n=j$.
Since we have 
$g_{j+1}(\cdisc_{j+1})\subset U_{j+1}$ and 
$g_{j+1}(\bar{\disc_{j+1}\setminus \Omega_j})\cap Y=\emptyset$ and since
by (V) for $n=j$ and by (c) we get that
$g_{j+1}(\bar{\Omega_{j}\setminus V})\cap Y=\emptyset$ it follows that
$g_{j+1}(\bar{\disc_{j+1}\setminus V})\subset U_{j+1}\setminus Y$.
This together with (II) for $n=j$ implies that there is
$\e_{j+1}$, $0<\e_{j+1}<\e_j$, so small that
(V) holds for $n=j+1$ and that (VI) holds for $n=j$.
The construction is finished.

\medskip

Denote by $\Omega$ the union
$\cup_{j=1}^\infty\Omega_j$. 
As $\Omega$ is a union of an increasing sequence of simply connected open sets
it is simply connected and therefore conformally equivalent to the unit disc.
It follows by (III) that for $\z\in\Omega$ the sequence $g_n(\z)$ is Cauchy
with respect to the complete metric $d$ therefore it converges to $g(\z)$. 
Since the convergence is uniform on compact sets it
follows that the map $g$ is \holo\ on $\Omega$. 

Next we show that the map $g$ and the domain $\Omega$ have all the required properties.
Fix $j\in\N\cup\{0\}$.
It follows by (III) that
\begin{eqnarray}
d(g(\z),g_{j+1}(\z))&\le& d(g_{j+1}(\z),g_{j+2}(\z))+
d(g_{j+2}(\z),g_{j+3}(\z))+\cdots<\nonumber\\
&<&\frac{\e_{j+1}}{2^{j+1}}+\frac{\e_{j+2}}{2^{j+2}}+\cdots<\e_{j+1}\ (\z\in \Omega_{j+1}).
\label{disk1}
\end{eqnarray}
Thus for $\z\in \Omega_{j+1}\setminus \Omega_j$
it holds by (VI) that $g(\z)\notin U_{j-1}$. This implies that 
$g\colon\Omega\to X$ is a proper map.
To prove that $g(\Omega\setminus V)$ avoids $Y$, choose $\z\in\Omega\setminus V$.
There is $j\in\N$ so large that $\z\in \Omega_{j+1}$. 
It follows by (\ref{disk1}) and by (V) that $g(\z)\notin Y$.
By (\ref{disk1}) for $j=0$ we obtain that $d(g(\z),f(\z))<\e$ for $(\z \in V)$.
We get by (IV) that $g(0)=f(0)$ and $g'(0)=f'(0)$.
This completes the proof.
\end{proof}

{\it Proof of Theorem \ref{thm1.2}.\,}
There are $r$ and $R$, $0<r<R<1$, \st\ $V\ss r\disc$.
One can choose $\e_0$, $0<\e_0<\e$, so small that
\begin{equation}
\text{for } z\in X \text{ such that } d(z,f(r\disc\setminus V))<\e_0
\text{ it holds that }z\notin Y.\label{d1}
\end{equation}
There is $\d>0$ so small that 
\begin{eqnarray}
V\ss(r-\d)\disc\subset (r+\d)\disc\ss R\disc ,\label{d2}\\
\z\in r\disc,\ \z'\in\disc \text{ such that }|\z-\z'|<\d\text{ then }
d(f(\z),f(\z'))<\tfrac{\e_0}{2}. \label{d3}
\end{eqnarray}
Let $\Omega_0=\emptyset$ and choose an increasing sequence $\{R_n\}$ of positive numbers
converging to $1$ with $R_1>R$.  
We shall construct inductively an increasing sequence of simply connected domains
$\{\Omega_n\}$ \st\ $R_n\disc\cup\Omega_{n-1}\ss\Omega_n\ss\disc$,
a decreasing sequence of positive numbers $\{\e_n\}$, $\e_1<\tfrac{\e_0}{2}$,
 and a sequence of proper holomorphic maps
$g_n\colon\Omega_n\to X$ \st\ 
\begin{enumerate}
\item[(a)] $g_n(\z)\notin Y$ for $\z\in\Omega_n\setminus V$,
\item[(b)] $d(g_n(\z),f(\z))<{\e_n}$ for $\z\in R_n\disc\cup\Omega_{n-1}$,
\item[(c)] $g_n(0)=f(0)$ and $g_n'(0)=f'(0)$.
\end{enumerate}
Assume that we have already constructed $\Omega_n$ and $\e_n$  $(0\le n\le j)$ and
$g_n$ $(1\le n\le j)$ for some
$j\in\N\cup \{0\}$.
One can choose $\e_{j+1}$, $0<\e_{j+1}<\tfrac{\e_j}{2}$, 
with the following property
\begin{equation}
\text{for } z\in X \text{ such that } d(z,f((R_{j+1}\disc\cup\Omega_{j})\setminus V))<\e_{j+1}
\text{ it holds that }z\notin Y .\label{dn1}
\end{equation}
Using Lemma \ref{lema} for $V=R_{j+1}\disc\cup\Omega_j$ and $\e=\e_{j+1}$
we obtain a simply connected domain $\Omega_{j+1}$,
$R_{j+1}\disc\cup\Omega_j\ss\Omega_{j+1}\ss\disc$, and a proper holomorphic
map $g_{j+1}\colon\Omega_{j+1}\to X$ which satisfy (b) and (c) for $n=j+1$ and
it holds that $g_{j+1}(\z)\notin Y$ for $\z\in\Omega_{j+1}\setminus 
(R_{j+1}\disc\cup\Omega_j)$.
By (\ref{dn1}) and (b) we get $g_{j+1}(\z)\notin Y$ for $\z\in (R_{j+1}\disc\cup\Omega_j)
\setminus V$, thus (a) holds for $n=j+1$. This completes the construction.

Note that $\cup_n\Omega_n=\disc$.
Caratheodory kernel theorem \cite{Car,Pom} implies that the sequence of
conformal maps $h_n\colon\disc\to\Omega_n$, \st\ $h_n(0)=0$, $h_n'(0)>0$,
converges uniformly on compact sets to identity.
Choose $n$ so big that 
\begin{equation}
|h_n(\z)-\z|<\d\ \ (\z\in r\disc).
\label{d4}
\end{equation}
Let $g=g_n\circ h_n$. By the above $g\colon\disc\to X$
is a proper holomorphic map and (c) implies (iii).
Take $\z\in r\disc$. By (\ref{d4}) and (\ref{d2}) we get that
$h_n(\z)\in R\disc$ which by (b) implies that 
$d(g_n(h_n(\z)),f(h_n(\z)))<\tfrac{\e_0}{2}$.
It follows by (\ref{d4}) and (\ref{d3}) that $d(f(h_n(\z)),f(\z))<\tfrac{\e_0}{2}$.
Therefore $d(g(\z),f(\z))<{\e_0}$ $(\z\in r\disc)$.
This proves (ii) and for $\z\in r\disc\setminus V$ this together with (\ref{d1}) implies that
$g(\z)\notin Y$.
Choose $\z\in\disc\setminus r\disc$.
By (\ref{d4}) it follows from Rouch\'{e}'s theorem that
$(r-\d)\disc\subset h_n(r\disc)$ and 
thus we get by (\ref{d2}) that $h_n(\z)\notin V$.
By (a) it follows that $g(\z)\notin Y$. This proves (i). The proof is complete.
{\hfill $\Box$ \\ \par}

{\it Acknowledgments.} The author wishes to thank 
J. Globevnik for many helpful
discussions while working on this paper. 
She would also like to thank 
F. Forstneri\v{c} whose suggestions led to considerable
simplification of the proof of Lemma \ref{lema}.
Thanks to M. \v{C}erne for pointing out Corollary \ref{cor1}.
The author is grateful to the referee for his/her accurate remarks and comments.

This research has been supported in part by 
the Ministry of Education, Science and Sport of Slovenia
through research program Analysis and Geometry, Contract No. P1-0291.

\medskip

\noindent Institute of Mathematics, Physics and Mechanics

\noindent University of Ljubljana

\noindent Jadranska 19

\noindent SI-1000 Ljubljana, Slovenia

\noindent e-mail address: Barbara.Drinovec@fmf.uni-lj.si


\begin{thebibliography}{999}

\bibitem[Ahl]{Ahl} L. Ahlfors, {\it Open Riemann surfaces and extremal problems on compact
subregions}, Comment. Math. Helvetici {\bf 24} (1950), 100--123.

\bibitem[Ale]{Ale} H. Alexander, {\it On a problem of Julia}, Duke Math. J. {\bf 42}
(1975), 327--332.

\bibitem[Car]{Car} C. Caratheodory, {\it Untersuchungen \"uber die konformen Abbildungen
von festenund ver\"anderlichen Gebieten},  {Math.\ Ann.} {\bf 72} (1912),  {107--144}.
 
\bibitem[Ch1]{Ch1} E. M. Chirka, {\it Rado's theorem for CR-mappings of hypersurfaces} (Russian), Mat. Sb. {\bf 185}
 (1994), 125--144, translation in Russian Acad. Sci. Sb. Math. {\bf 82} (1995), 243--259.
 
\bibitem[Ch2]{Ch2} E. M. Chirka, {\it Complex analytic sets} (in English), Kluwer, Dordrecht, 1989.

\bibitem[Col]{Col} M. Col\c{t}oiu, {\it Complete locally pluripolar sets}, 
J. Reine Angew. Math. {\bf 412} (1990), 108-112.

\bibitem [Dri]{Dri} B. Drinovec Drnov\v sek, {\it Proper holomorphic discs avoiding closed convex sets},
{Math.\ Z.} {\bf 241} (2002),  {593--596}.


\bibitem[FG1]{FG1} F. Forstneri\v{c}, J. Globevnik, {\it Discs in pseudoconvex
domains}, Comment. Math. Helvetici {\bf 67} (1992), 129--145.

\bibitem[FG2]{FG} F. Forstneri\v{c}, J. Globevnik, {\it Proper holomorphic discs in 
$\ctwo$}, Math. Res. Lett. {\bf 8} (2001), 257--274.

\bibitem[Glo]{Glo} J. Globevnik, {\it Discs in Stein manifolds}, Indiana Univ. Math. J.
{\bf 49} (2000), 553--574.

\bibitem[H\"or]{Hor} L. H\"ormander, {\it An introduction to complex 
analysis in several variables}, Van Nostrand, Princeton, 1973.

\bibitem[Jul]{Jul} G. Julia, {\it Sur le domaine d'existence d'une fonction implicite d$\acute e$finie
par une relation enti$\grave e$re $G(x,y)=0$}, Bull. Soc. Math. France 
{\bf 54} (1926), 26--37.

\bibitem[Pom]{Pom} C. Pommerenke, {\it Univalent functions}, (Math. Lehrb., Bd. 25),
Vanderhoeck and Rupprecht, G${\rm \ddot{o}}$ttingen, 1975

\bibitem[Re1]{Re1} R. Remmert, {\it Projektionen analytischer Mengen}, 
Math. Ann. {\bf 130} (1956), 410--441.

\bibitem[Re2]{Re2} R. Remmert, {\it Holomorphe und meromorphe Abbildungen komplexer Räume},
Math. Ann. {\bf 133} (1957), 328--370.

\bibitem[Tsu]{Tsu} M. Tsuji, {\it Theory of meromorphic functions in a neighborhood of a closed
set of capacity zero}, Japanese J. Math. {\bf 19} (1944), 139--154.

\end{thebibliography}
\end{document}